\documentclass[12pt]{article}

\newif\ifrestated
\restatedfalse

\usepackage{booktabs} % for much better looking tables
\usepackage[table]{xcolor}
\usepackage{pifont}
\usepackage{soul}

\usepackage[utf8]{inputenc}
\usepackage[english]{babel}

\usepackage{geometry}
\usepackage{amsmath}
\usepackage{amssymb}
\usepackage{mathtools}
\usepackage{xfrac}
\usepackage{dsfont}

\usepackage{hyperref}
\usepackage{pdfpages}
\hypersetup{colorlinks=true,linkcolor=blue}

\usepackage[thmmarks,thref,hyperref]{ntheorem}
\usepackage{thmtools, thm-restate}

\theoremstyle{plain}
\theoremseparator{. }
\newtheorem{lemma}{Lemma}

\usepackage{pgf,tikz}
\usetikzlibrary{arrows,automata, patterns, calc}

\usepackage[hypcap=true]{caption}

\def\thedate{\today}
\def\me{I. Fabrici, J. Harant, T. Madaras, S. Mohr, R. Soták, \\C. T. Zamfirescu}

\def\departA{Ilmenau University of Technology, Department of Mathematics, Ilmenau, Germany}

\def\departB{Pavol Jozef Šafárik University, Institute of Mathematics, Košice, Slovakia}

\def\departC{Ghent University, Department of Applied Mathematics, Computer Science and Statistics, Ghent, Belgium}

\def\departD{Babe\c{s}-Bolyai University, Department of Mathematics, Cluj-Napoca, Roumania}

\def\titlevar{Long Cycles and Spanning Subgraphs of\\[1mm] Locally Maximal 1-planar Graphs}

\title\titlevar
\author\me
\date\thedate

\hypersetup{pdftitle=\titlevar,pdfsubject=Note,pdfauthor=\me}

\parskip\medskipamount %Abstand zwischen 2 Absätzen

\makeatletter
\newcommand\need[1]{\par \penalty-100 \begingroup % preserve \dimen@
   \dimen@\pagegoal \advance\dimen@-\pagetotal % space left
   \ifdim #1>\dimen@ % not enough space left
      \ifdim\dimen@>\z@ \vskip -\pagedepth plus 1fil \fi
      \break
   \fi \endgroup}

\let\circ\relax
\DeclareMathOperator\circ{circ}

\let\oldoverline\overline
\renewcommand{\overline}[2][3]{{}\mkern#1mu\oldoverline{\mkern-#1mu#2}}

\makeatother

\begin{document}
\begin{center}
{\bf \Large \titlevar }\\[3mm]
{\bf I.~Fabrici\footnote{\label{DAAD}Partially supported by DAAD, Germany (as part of BMBF) and the Ministry of Education, Science, Research and Sport of the Slovak Republic within the project 57447800.}\textsuperscript{\textnormal{,}}\footnote{\label{APVV}Partially supported by the Slovak Research and Development Agency under contract No.\ APVV-15-0116 and by the Science Grant Agency -- project VEGA 1/0368/16.}\textsuperscript{\textnormal{,a}},
J.~Harant\textsuperscript{\textnormal{\ref{DAAD},b}},
T.~Madaras\textsuperscript{\textnormal{\ref{DAAD},\ref{APVV},a}},
S.~Mohr\textsuperscript{\textnormal{\ref{DAAD},}}\footnote{\label{DFG}Gefördert durch die Deutsche Forschungsgemeinschaft (DFG) -- 327533333; \\partially supported by the grant 327533333 from the Deutsche Forschungsgemeinschaft (DFG, German Research Foundation).}\textsuperscript{\textnormal{,b}},
R.~Soták\textsuperscript{\textnormal{\ref{DAAD},\ref{APVV},a}},
C.~T.~Zamfirescu\footnote{Supported by a Postdoctoral Fellowship of the Research Foundation Flanders (FWO).}\textsuperscript{\textnormal{,c,d}}}\\
\textsuperscript{a} \departB \\
\textsuperscript{b} \departA \\
\textsuperscript{c} \departC \\
\textsuperscript{d} \departD \\
\end{center}

\hrulefill

\begin{abstract}
A graph is {\em $1$-planar} if it has a drawing in the plane such that
each edge is crossed at most once by another edge. Moreover, if this drawing has the additional property that for each crossing of two edges
the  end vertices of these edges induce a complete subgraph, then  the graph is {\em locally maximal $1$-planar}.
For  a $3$-connected locally maximal $1$-planar graph $G$, we show the existence of a spanning $3$-connected planar subgraph and prove that $G$
is hamiltonian if $G$ has at most three $3$-vertex-cuts, and that $G$ is traceable if $G$ has at most four $3$-vertex-cuts. Moreover, infinitely many non-traceable $5$-connected $1$-planar graphs are presented.
\end{abstract}

\bigskip


\section{Introduction and Results}

We use standard terminology of graph theory and consider  finite and simple graphs, where $V(G)$ and $E(G)$ denote the vertex set  and the edge
set of a graph $G$, respectively.
These graphs are represented by {\em drawings} in the plane, such that vertices are distinct points and edges are arcs (non-self-intersecting
continuous curves, i.\,e.\ open Jordan curves)
that join two points corresponding to their incident vertices.
The arcs are supposed to  contain no vertex points in their interior.
Such a drawing of a graph $G$ in the plane is denoted by $D(G)$. For more details on drawings of graphs in the plane, see~\cite{hudak2012properties,pach1997graphs}. If two edges  of $D(G)$ have an internal point in common, then these edges {\em cross} and we call the
pair of these edges a {\em crossing}, and the aforementioned internal point their {\em crossing point}.
It is easy to see that a drawing can be changed locally to  a different drawing with  fewer crossings
if two edges with a shared end vertex cross or
if two edges cross several times.    Thus, in the sequel we will consider  drawings  with the property that if two edges cross, then they do so exactly once
and their four end vertices are mutually distinct.

A graph $G$ is {\em planar} if there exists a drawing $D(G)$ of $G$ without crossings.
There are several different approaches generalizing the concept of planarity.
One of them is to allow a given constant number of crossings for each edge in a drawing $D(G)$.
In particular, if there exists a drawing $D(G)$ of a graph $G$ such that each edge of $D(G)$ is crossed at most once by another edge, then $G$ is
 {\em $1$-planar}. In this case we call $D(G)$ a {\em $1$-planar drawing} of  $G$.
This class of graphs was introduced by  Ringel~\cite{ringel1965sechsfarbenproblem} in connection with the simultaneous vertex-face coloring of plane graphs.
Properties of $1$-planar graphs are studied in \cite{fabrici2007structure,hudak2012properties,Korzhik,KorzhikMohar,pach1997graphs}.

Pach and Tóth~\cite{pach1997graphs} proved
that each $1$-planar graph on $n$ vertices, $n\ge 3$, has at most $4n - 8$ edges and this bound
is attained for every $n\geq 12$. As a consequence, a $1$-planar graph has a vertex of degree at most $7$, hence, it is at most $7$-connected. A
$1$-planar graph on $n$ vertices is {\em optimal} if it has exactly $4n - 8$ edges.
A  graph $G$ from a family $\mathcal{G}$ of graphs is {\em maximal} if $G + uv \notin \mathcal{G}$ for any
two non-adjacent vertices $u, v \in V(G)$.
It is remarkable that there exist maximal $1$-planar graphs on $n$ vertices that have significantly fewer than $4n-8$ edges~\cite{Brandenburg}.
Thus, in contrast to the planar case, the properties ``optimal'' and ``maximal'' are not the same for $1$-planar graphs. Obviously, an optimal
$1$-planar graph is also maximal $1$-planar. It is clear that a maximal planar graph is not necessarily maximal $1$-planar.

The length (number of vertices) of a longest cycle of a graph $G$ (also called \emph{circumference} of $G$) is denoted by $\circ(G)$.
If $\circ(G)=n$ for a graph $G$ on $n$ vertices, then $G$ is \emph{hamiltonian} and a longest cycle of $G$ is a {\em hamiltonian
cycle}. In the same vein, a graph is \emph{traceable} if it contains a path visiting every vertex of the graph.

In \cite{hudak2012properties}, it is proved that an optimal $1$-planar graph is hamiltonian. This is in sharp contrast with the family of planar
graphs since Moon and Moser \cite{MoonMoser} constructed infinitely many maximal planar graphs $G$ with $\circ(G)\le 9 |V(G)|^{\log
_32}$ (in fact, Moon and Moser even showed that every path in $G$ is strictly shorter than the aforementioned length). 
It is
known that a maximal planar graph on $n\ge 4$ vertices is $3$-connected. In~\cite{hudak2012properties}, maximal $1$-planar graphs with vertices of degree $2$ are constructed and
it remained open there whether every $3$-connected maximal $1$-planar graph is hamiltonian. Moreover, the question arises
whether such a construction as the one of Moon and Moser is also possible in the class of $3$-connected maximal $1$-planar graphs.
An answer to both questions is given by  \thref{W1}  which has the consequence that there are positive constants $c$ and $\alpha \le\log_32<1$  such that infinitely many $3$-connected maximal $1$-planar graphs $G$ with $\circ(G)\le c\cdot |V(G)|^\alpha$  exist.

\begin{restatable}{theorem}{thmone}
\label{W1}
If $H$ is a maximal planar graph on $n\ge 4$ vertices, then there is a $3$-connected maximal $1$-planar graph $G$ on $7n-12$ vertices such that
$\circ(G)\le 4\cdot \circ(H)$.
\end{restatable}

For an arbitrary $1$-planar drawing $D(G)$ of a graph $G$, let $D^\times(G)$ be the plane graph obtained from $D(G)$ by turning all crossings into new $4$-valent vertices. If $uv$ and $xy$ are two crossing edges of $D(G)$, then let $c$ be the
vertex of $D^\times(G)$ corresponding to the crossing point of $uv$ and $xy$. Let $\alpha$ be the face of $D^\times(G)$ such that $ucx$
is a subpath of the facial walk of $\alpha$ in $D^\times(G)$. If $ux$ is an edge of $G$ and $ux$ is crossed by another edge in $D(G)$, then it is possible to redraw the edge $ux$ in $D(G)$ such that $ux$ lies in the region of $D(G)$ corresponding to the face $\alpha$ of $D^\times(G)$. It follows that $ux$ is not crossed by another edge in $D(G)$ anymore. Thus, in the sequel we will consider $1$-planar drawings $D(G)$ of a graph $G$ with  the property that if $uv$ and $xy$ are crossing edges of $D(G)$, then the edge
$xu$ (if it exists) is not crossed by another edge in $D(G)$.

Now we will consider much wider classes of $1$-planar graphs than the class of maximal $1$-planar graphs.
Let $K_4^-$ be the graph obtained from the complete graph $K_4$ on four vertices by removing one edge.
Given a $1$-planar drawing $D(G)$ of a  graph $G$, we call a crossing of $D(G)$ \emph{full} or \emph{almost full} if the four end vertices of its edges induce a $K_4$ or a $K_4^-$, respectively.

If for a $1$-planar graph $G$ there exists a $1$-planar drawing $D(G)$ such that all crossings of $D(G)$ are full or almost full, or all crossings of $D(G)$ are full, then, in the first case, we call $G$ {\em weakly locally maximal $1$-planar} and $D(G)$ a  {\em weakly locally maximal $1$-planar drawing}
of $G$ or, in the second case, $G$ {\em  locally maximal $1$-planar} and $D(G)$ a  {\em  locally maximal $1$-planar drawing}
of $G$, respectively. Obviously, a planar graph is locally maximal $1$-planar and it can easily be seen that a maximal $1$-planar graph is also locally maximal
$1$-planar and that a locally maximal $1$-planar graph is also weakly locally maximal
$1$-planar.
For a positive integer $k\ge 2$, Figure~\ref{lm}  shows a graph on $4k$ vertices which is  locally maximal $1$-planar, obviously not  maximal $1$-planar, and also not planar since it contains a subdivision of $K_5$ with major ($4$-valent) vertices $u_1,x_1,y_1,z_1,x_k$.

\begin{figure}[h]
\begin{center}
\begin{tikzpicture}[scale=0.2]
\foreach \w in {20,40,...,360} {

\fill (\w:7cm) circle (0.3cm);

\fill (\w:10cm) circle (0.3cm);

\draw[thick] (\w:7cm) -- (\w+20:7cm);
\draw[thick] (\w:10cm) -- (\w+20:10cm);
\draw[thick] (\w:7cm) -- (\w:10cm);

}

\foreach \w in {40,80,...,360} {

\draw[thick] (\w:7cm) -- (\w+20:10cm);
\draw[thick] (\w:10cm) -- (\w+20:7cm);

}

\draw (60:5.5cm) node {$x_k$};
\draw (60:11.5cm) node {$y_k$};
\draw (40:5.5cm) node {$z_k$};
\draw (40:11.5cm) node {$u_k$};

\draw (100:5.5cm) node {$x_1$};
\draw (100:11.5cm) node {$y_1$};k
\draw (80:5.5cm) node {$z_1$};
\draw (80:11.5cm) node {$u_1$};

\draw (140:5.5cm) node {$x_2$};
\draw (140:11.4cm) node {$y_2$};
\draw (120:5.5cm) node {$z_2$};
\draw (120:11.5cm) node {$u_2$};

\end{tikzpicture}
\caption{}\label{lm}
\end{center}
\end{figure}

Whitney~\cite{whitney1931theorem} showed that a $4$-connected maximal planar graph is hamiltonian. Later Tutte~\cite{tutte1956theorem}
proved that an arbitrary $4$-connected planar graph has a hamiltonian cycle. We
remark that  non-hamiltonian $4$-connected $1$-planar graphs are constructed in \cite{hudak2012properties}. In order to formulate the next result, we recall that for an infinite family ${\cal G}$ of graphs, its \emph{shortness exponent} is defined as $$\sigma({\cal G})=\liminf \limits_{G\in {\cal G}}\frac{\log \circ(G)}{\log |V(G)|}.$$ See \cite{GW} for details concerning the theory of shortness exponents. We are able to prove the following theorem---however, it remains open whether a non-hamiltonian $6$-connected $1$-planar graph exists.

\begin{restatable}{theorem}{thmtwo}
\label{thm:5conn}
There are infinitely many non-traceable $5$-connected weakly locally maximal $1$-planar graphs. Moreover, for the class $\Gamma$ of $5$-connected $1$-planar graphs we have $\sigma(\Gamma)\le \frac{\log 20}{\log 21}$.
\end{restatable}

We can infer from \thref{thm:5conn} that for an arbitrary $\varepsilon >0$ there is a $5$-connected $1$-planar graph $G$ such that $\circ(G)<\varepsilon\cdot|V(G)|$.

It is also not known whether every $7$-connected $1$-planar graph is
hamiltonian (see \cite{hudak2012properties}), so the intriguing question whether an analog of Tutte's theorem holds for the family of $1$-planar graphs remains open.

By \thref{W1}, $3$-connected maximal $1$-planar graphs are far away from being hamiltonian in general---nonetheless \thref{thm:3cut} and \thref{thm:almost} both imply that
a $4$-connected locally maximal $1$-planar graph is hamiltonian, i.\,e.\ Whitney's theorem  can  be  extended to the class of $4$-connected locally maximal
$1$-planar graphs.
For an overview of the minimum sufficient connectivity that leads to hamiltonicity for the different kinds of 1-planar maximality discussed in this article, we refer the reader to Table~\ref{tab:overview} at the end of this paper.

As an extension of Tutte's theorem, it is proved in \cite{bz} that a $3$-connected planar graph with at most three $3$-cuts is hamiltonian. (In this paper, all cuts are vertex-cuts.) By \thref{thm:5conn}, this result cannot be extended to the class of $3$-connected weakly locally maximal $1$-planar graphs, however, \thref{thm:3cut} shows that the assertion is true for $3$-connected locally maximal $1$-planar graphs.

\begin{restatable}{theorem}{thmthree}
\label{thm:3cut}
A $3$-connected locally maximal $1$-planar graph with at most three $3$-cuts is hamiltonian. Furthermore, every $3$-connected locally maximal $1$-planar graph with at most four $3$-cuts is traceable.
\end{restatable}

By \thref{thm:5conn}, there are infinitely many non-hamiltonian $5$-connected weakly locally maximal $1$-planar graphs. \thref{thm:almost} shows that the situation changes if the number of almost full crossings in a  weakly locally maximal $1$-planar drawing of a graph is not too large, even if this graph is only $4$-connected.

\begin{restatable}{theorem}{thmfour}
\label{thm:almost}
If a $4$-connected graph has a weakly locally maximal $1$-planar drawing with at most three almost full crossings, then it is hamiltonian. Moreover, if a $4$-connected graph has a weakly locally maximal $1$-planar drawing with at most four almost full crossings, then it is traceable.
\end{restatable}

Chen and Yu~\cite{ChenYu} showed that there is a constant $c$ such that $\circ(G)\ge c\cdot |V(G)|^{\log _32}$ for an arbitrary
$3$-connected planar graph $G$. 
By~\thref{Subgraph}, we
show that the  extension of the result of Chen and Yu (and of any other  result concerning lower bounds on the length of a longest cycle of a
$3$-connected planar graph) to $3$-connected locally maximal $1$-planar graphs is possible. Moreover, by~\thref{Subgraph}, every result on
the existence of a certain subgraph of a  $3$-connected planar graph is also true for a $3$-connected locally maximal $1$-planar graph. Examples
are the  results of Barnette~\cite{Barnette} that a $3$-connected planar graph has a spanning tree of maximum degree at most $3$, and of Gao~\cite{Gao} that a $3$-connected planar graph has a spanning $2$-connected subgraph of maximum degree at most~$6$.

\begin{restatable}{theorem}{thmfive}
\label{Subgraph}
Each  $3$-connected locally maximal $1$-planar graph has a  $3$-connected planar spanning subgraph.
\end{restatable}

In Figure~\ref{lm}, a $4$-connected locally maximal $1$-planar graph is presented. Because it is non-planar and $4$-regular, it cannot contain a 
$4$-connected planar spanning subgraph. Thus, \thref{Subgraph} is best possible in this sense.

One obtains a weakly locally maximal $1$-planar graph $G$ (and one of its weakly locally maximal $1$-planar drawings) if the edges $x_1y_1,\dots, x_ky_k$ are removed from the graph of Figure~\ref{lm}. Assume this graph $G$ contains a  $3$-connected planar spanning subgraph $H$. Since $H$ has minimum degree at least three, all edges incident with a vertex from $\{x_1,\dots,x_k,y_1,\dots,y_k\}$ belong to $H$. Thus, the graph $H'$ obtained from $H$ by removing the edges $u_1z_1,\dots,u_kz_k$ is a subgraph of $H$. If $k\ge 3$, then it is easy to see that $H'$ contains a subdivision of $K_{3,3}$ with major ($3$-valent) vertices  $x_1,y_1,x_k$ and $u_1,z_1,x_2$, a contradiction to the planarity of $H$. It follows that \thref{Subgraph} is also best possible in the sense that ``locally maximal $1$-planar'' cannot be replaced with ``weakly locally maximal $1$-planar''.

\section{Proofs}\label{sec:known}

For the proof of \thref{W1}, we need a result of Bachmaier et al.~\cite{Brandenburg1}. 
Let $D(G)$ be an arbitrary $1$-planar drawing of a $1$-planar graph $G$. 
Two subgraphs  $H_1$ and $H_2$ of $G$ are said to be {\em $k$-sharing} if $H_1$ and $H_2$ have at least $k$ vertices in common. 
Moreover, $H_1$ and $H_2$  {\em share a crossing} in $D(G)$ if there are edges $e_1\in E(H_1)$ and $e_2\in E(H_2)$ that cross in $D(G)$. 
The following lemma holds:

\begin{lemma}[Lemma~7 from~\cite{Brandenburg1}]\label{lem:k5share}
Let $D(G)$ be a $1$-planar drawing of a $1$-planar graph $G$. 
If two subgraphs of $G$ both isomorphic to $K_5$ share a crossing in $D(G)$, then they are $3$-sharing. 
\end{lemma}

\bigskip

{\bf Proof of}
\ifrestated
\vspace{-1.6ex}

\thmone*
\else
{\bf\thref{W1}. }
\fi

We construct $G$ from $H$ such that $H$ is a subgraph of $G$. Therefore,  the vertices  of $H$ are said to be {\em old} and these in
$V(G)\setminus V(H)$  to be {\em new}. \\
It is well-known that a simple maximal planar graph on at least $4$ vertices is $3$-connected. Whitney~\cite{whitney1932theorem} (see also \cite{Fleischner}) proved, that a $3$-connected planar graph has a unique (up to the choice of the outer face) planar drawing. \\
Let $D_0(H)$ be (in this sense) the unique planar drawing of $H$. Figure \ref{K6} shows a $1$-planar embedding of $K_6$. A drawing $D_0(G)$ of $G$ is obtained from $D_0(H)$ by inserting into each face of $D_0(H)$ with  boundary vertices $u,v,$ and $w$
a triangle with three new vertices $a,b,$ and $c$,  and completed by further nine edges as shown in Figure~\ref{K6}.

\begin{figure}[h]
\begin{center}
\begin{tikzpicture}[scale=0.18]
\foreach \w in {90,210,330} {

\fill (\w:4cm) circle (0.5cm);
\fill (\w:10cm) circle (0.5cm);

\draw[thick] (\w:4cm) -- (\w+120:4cm);
\draw[thick] (\w:10cm) -- (\w+120:10cm);
\draw[thick] (\w:4cm) -- (\w+120:10cm);
\draw[thick] (\w:4cm) -- (\w-120:10cm);
\draw[thick] (\w:4cm) -- (\w:10cm);

}
\draw (90:2cm) node { $a$};
\draw (210:2cm) node { $b$};
\draw (330:2cm) node { $c$};
\draw (90:12cm) node { $u$};
\draw (210:12cm) node { $v$};
\draw (330:12cm) node { $w$};

\end{tikzpicture}
\caption{}\label{K6}
\end{center}
\end{figure}

Obviously, $D_0(G)$ is a
   $1$-planar drawing of $G$, hence, $G$ is $1$-planar.\\
As $H$ is maximal planar with $2n-4$ faces,  $G$ has $n+3(2n-4)=7n-12$ vertices.\\
Let $C_G$ be a longest cycle of $G$. If $C_G$ has at most eight vertices, then $\circ(G)=|V(C_G)|\le 4\cdot \circ(H)$ is true.
If $P$ is a subpath of $C_G$ connecting two old vertices $u$ and $v$ such that $V(P)\setminus \{u,v\}$
does not contain an old vertex, then $u$ and $v$ are vertices of a face of $H$  and  $|V(P)\setminus \{u,v\}|\le 3$ (see Figure \ref{K6}). In case $V(P)\setminus \{u,v\}\neq \emptyset$, let $Q$ be the $u$ and $v$ connecting path obtained from $C_G$ by removing $V(P)\setminus \{u,v\}$. Note that $Q$ contains at least three old vertices because $|V(C_G)|\ge 9$. We add the edge $uv$ to $Q$ and the resulting graph is again a cycle of $G$ containing fewer new vertices than $C_G$.
Repeating this step, we obtain a
cycle $C_H$ of $H$ and it follows that $\circ(G)=|V(C_G)|\le 4\cdot |V(C_H)|\le 4\cdot \circ(H)$.

Now we show that $G$ is $3$-connected. 
Therefore,
consider a minimal cut $S$ of $G$. Since the neighborhood of each new vertex forms a complete graph, $S$ does not contain new vertices, hence,
$S\subset V(H)$. If $S$ is also a cut of $H$, then $|S|\ge 3$ since a simple maximal planar graph is $3$-connected.
If $|S|<3$, then, since
$H-S$ is still connected, $G-S$ has a component consisting  of new vertices only. This is impossible since each new vertex has three old neighbors.

For the proof of \thref{W1}, it remains to show that $G$ is maximal $1$-planar.
Therefore, let $D(G)$ be an arbitrary $1$-planar drawing of $G$. We will prove  {\em (i)}. \\

{\em (i) An edge of $H$ is not crossed in $D(G)$.}\\

Let $uv$ be an  edge of $H$ (see Figure~\ref{K6}) and $\{a,b,c\}$ and $\{a',b',c'\}$ be the disjoint sets of new vertices being inserted into the two faces of $H$ both incident with $uv$, respectively. The subgraphs $G[\{u,v,a,b,c\}]$ and $G[\{u,v,a',b',c'\}]$ of $G$ both are isomorphic to $K_5$. Assume there is an edge $e$ of $G$ that crosses $uv$ in $D(G)$. Since each edge of $G$ is an edge of a subgraph isomorphic to $K_5$, let $K(e)$ be a $K_5$-subgraph of $G$ containing $e$.\\
By \thref{lem:k5share}, it follows that $G[\{u,v,a,b,c\}]$ and $K(e)$ and also $G[\{u,v,a',b',c'\}]$ and $K(e)$ are $3$-sharing. This implies $V(K(e))\cap \{a,b,c\}\neq \emptyset$ and $V(K(e))\cap \{a',b',c'\}\neq \emptyset$.  Since there is no edge between the sets $\{a,b,c\}$ and $\{a',b',c'\}$, this contradicts the completeness of $K(e)$, and {\em (i)} is proved.\\

By {\em (i)}, the restriction $D(H)$ of $D(G)$ to $V(H)$ is a planar embedding of $H$. The planar embedding of $H$ is unique, thus, $D(H)=D_0(H)$. Consider a face $F$ of $D_0(H)$ with boundary vertices $u,v$, and $w$. Since $H$  has at least four vertices and again by {\em (i)}, the three new vertices $a,b,$ and $c$ all adjacent to $u,v,$ and $w$ lie in the interior of the face $F$ and, up to permutation of $a,b,c$, the situation of  Figure \ref{K6} occurs. Thus, we may assume {\em (ii)}. \\

{\em (ii) $G$ has the unique $1$-planar embedding $D_0(G)$.}\\

To show that $G$ is maximal $1$-planar, assume to the contrary that there are nonadjacent vertices $x$ and $y$ such that the graph $G+xy$ obtained from $G$ by adding the edge $xy$ is $1$-planar. Therefore, let $D(G+xy)$ be a $1$-planar drawing of $G+xy$. If $xy$ is removed from $D(G+xy)$, then we obtain a $1$-planar embedding of $G$ and this embedding is $D_0(G)$ because of {\em (ii)}. Thus, we may assume that in $D_0(G)$ the edge $xy$ can be added in such a way that the resulting drawing is still $1$-planar.

If $x$ is new, then let $x=a$ (see Figure \ref{K6}). Since $xy\notin E(G)$, $y$ is not among the six vertices in Figure \ref{K6} and the edge $xy$ has to cross at least one of the edges  $ub,uc,bw,cv$ in $D(G+xy)$, but each of them is already crossed in $D_0(G)$, a contradiction.

If  $x$ and $y$ are old vertices, then, because $H$ is maximal planar,  $H+xy$ is not planar anymore. Thus, $xy$ crosses an edge $e\in E(H)$ in $D(G+xy)$. Let $e=uv$ (see Figure \ref{K6}), then  $xy$ crosses  the edge $av$ or $bu$ in $D(G+xy)$. However, $av$ and $bu$  cross each other in $D_0(G)$, again a contradiction,
and \thref{W1} is proved.
\hfill $\blacksquare$

\bigskip

{\bf Proof of}
\ifrestated
\vspace{-1.6ex}

\thmtwo*
\else
{\bf\thref{thm:5conn}. }
\fi

\begin{figure}[h]
\begin{center}
\includegraphics[scale=.6]{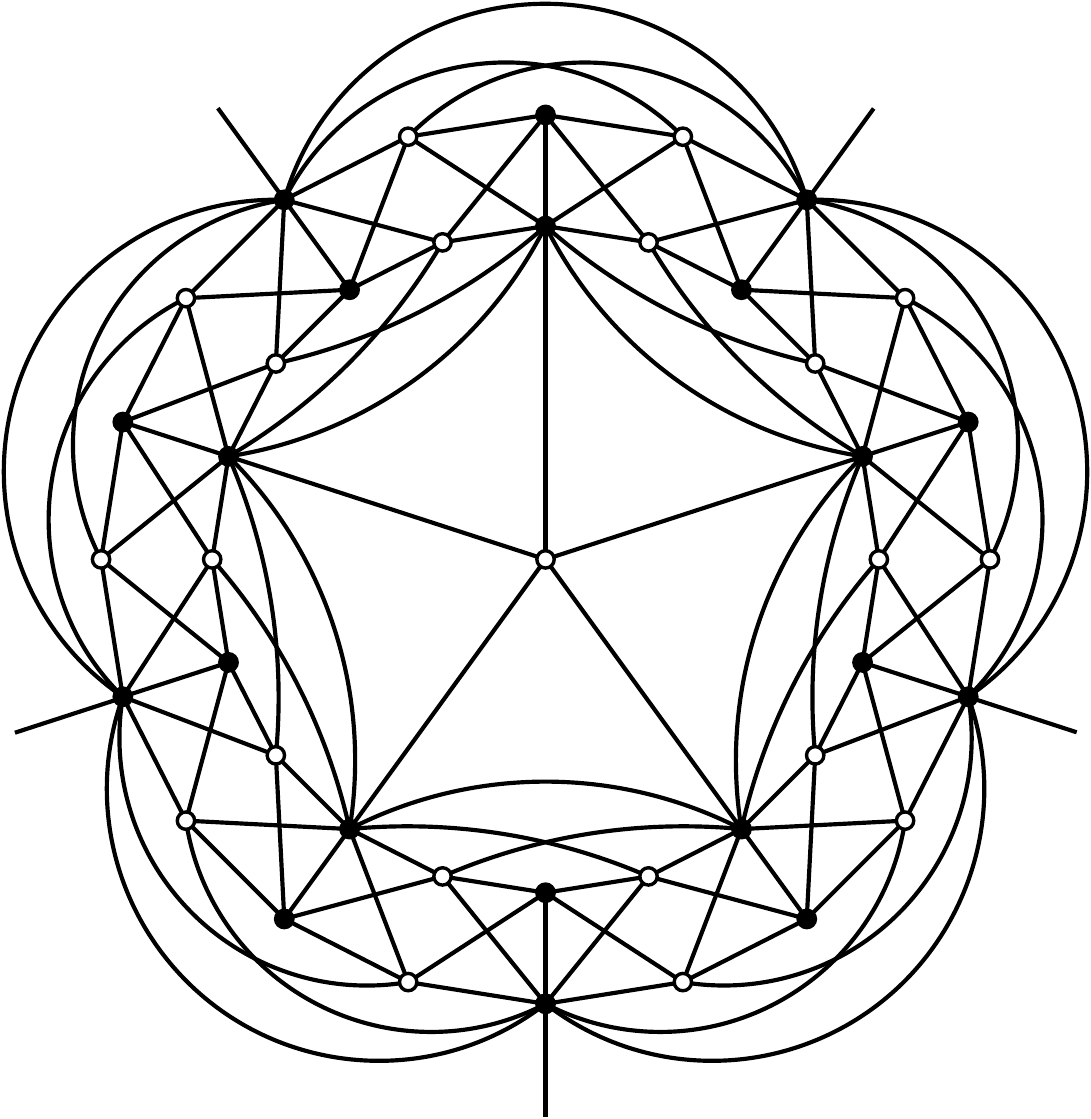}

\caption{The structure $H$}\label{5conn}
\end{center}
\end{figure}

Consider the structure $H$ shown in Figure~\ref{5conn}, add a new white vertex, and join the five half-edges of $H$ to this new vertex. We obtain a weakly locally maximal $1$-planar graph $G$. It is not difficult to see that $G$ is $5$-connected. Moreover, since $G$ contains $20$ black vertices and $22$ white vertices, and the set of white vertices is independent, it follows that $G$ is non-traceable. This construction can be generalized easily to obtain a weakly locally maximal $1$-planar graph  containing $4k$ black vertices and $4k+2$ independent white vertices, where $k$ is an arbitrary integer at least five.\\
Now construct the 5-connected weakly locally maximal $1$-planar graph $G_0$ from the graph $H$ by removing the five half-edges. Starting with $G_0$, we construct an infinite sequence $\{G_i\}$ for $i\ge 0$ of $5$-connected $1$-planar graphs as follows.
 Let $G_i$ be already constructed and $G_{i+1}$ be obtained by replacing each white vertex $v$ of $G_i$ with a copy $H_v$ of $H$ (Figure~\ref{5conn}) and connecting the five half-edges of $H$ with the five neighbors of $v$ in $G_i$. Let $M_i=\{H_v~|~v \text{~is white in~} G_{i-1}\}$ and $w_i$ be the number of white vertices of $G_i$.
Then $|M_{i+1}|=w_{i}$, $w_0=21$, and $w_{i+1}=21\cdot w_i$, thus, $|V(G_i)|>w_i=21^{i+1}$.

 Let $T_i$ be a longest closed trail of $G_i$ visiting each black vertex of $G_i$ at most once and let $t_i=|V(T_i)|$. Note that $T_i$ visits a white vertex $v$ of $G_i$ at most twice, since $v$ has degree five in $G_i$. Since a longest cycle of $G_i$ is also a closed trail of $G_i$, it follows that $\circ(G_i)\le t_i$ for all $i\ge 0$.

For $i\ge 1$, let
 $h_i=|\{H\in M_i~|~V(H)\cap V(T_i)\neq \emptyset\}|$ be the number of copies of $H$ in $G_i$ visited by $T_i$ at least once; it is easy to see that $h_i\ge 2$.
 Moreover, since the $21$ white  vertices of $H\in M_i$ are independent and a half-edge of $H$ is incident with a black vertex of $H$, it follows that  $V(T_i)\cap V(H)$ contains at most $19$ of the $21$ white vertices of $H$. Thus, $|V(T_i)\cap V(H)|\le 39$ and
 $t_i\le b_i + 39h_i$, where $b_i$ denotes the number of vertices of $T_i$ not belonging to some $H\in M_i$. Let $T$ be the closed trail of $G_{i-1}$ obtained from $T_i$ by shrinking all $H\in M_i$ to white vertices of $G_{i-1}$ again.

  Then $t_{i-1}\ge |V(T)|= b_i+h_i$  leads to $t_i-20t_{i-1}\le -19b_i+19h_i$ and, because all  $H\in M_i$ have distance at least $2$ in $G_i$,  $b_i\ge h_i$. Therefore, $t_i\le 20t_{i-1}$, which gives $\circ(G_i)\le t_i\le 20^it_0$. Finally, since
 $$\frac{\log \circ(G_i)}{\log |V(G_i)|}\le \frac{\log t_i}{\log |V(G_i)|}<\frac{\log 20+\frac{1}{i}\log t_0}{\log 21+\frac{1}{i}\log 21} \quad {\rm for} \quad i \ge 1,$$ we have
 $$\sigma(\Gamma)\le \lim\limits_{i \rightarrow \infty}\frac{\log \circ(G_i)}{\log |V(G_i)|}\le  \frac{\log 20}{\log 21}.$$
  \hfill $\blacksquare$

\bigskip

For the proofs of  \thref{thm:3cut} and \thref{thm:almost} we need the three forthcoming lemmas.

\begin{lemma}\label{lem:struct}
Let $t$ be a non-negative integer and $G$ be a non-planar $3$-connected weakly locally maximal $1$-planar graph that has a  weakly locally maximal
$1$-planar drawing with $t$ almost full crossings.
Furthermore, among all weakly locally maximal $1$-planar drawings of $G$ with at most $t$ almost full crossings let $D(G)$ be chosen with minimum
number of crossings.
Let $G'$ with a drawing $D(G')$ be constructed by turning an arbitrary crossing $X$ of $D(G)$ into a new $4$-valent vertex $v$.
\\
Then \\
(i)
$G'$ is weakly locally maximal $1$-planar and $D(G')$ has one crossing less than $D(G)$. \\
(ii) 
If $X$ is almost full, then $G'$ has a  weakly locally maximal $1$-planar drawing with at most $t-1$ almost full
crossings.
Otherwise, $G'$ has a  weakly locally maximal $1$-planar drawing with at most $t$ almost full crossings.
\\
(iii)  $G'$ is $3$-connected.
\\
(iv) Let $S$ be a $3$-cut of $G'$. 
If $S\subseteq V(G)$, then $S$ is also a $3$-cut of $G$.
If $v\in S$, then $X$ is almost full, the neighborhood of $v$ in $G'$ forms a path on vertices $a,b,c,d$ that appear in this order, and there is $z\in V(G)\setminus N_{G'}(v)$ such that $G$ has a $3$-cut $S'=\{b,c,z\}$  which separates $a$ and $d$. 
\end{lemma}

{\bf Proof of \thref{lem:struct}.}\\
Obviously, $G'$ is weakly locally maximal $1$-planar
(remember that 
all considered $1$-planar drawings $D(G)$ of a graph $G$ have the property that 
if two edges $uv$ and $xy$ are crossing edges of $D(G)$, then the edge $xu$ (if it exists) is not crossed by 
another edge in $D(G)$). 
Furthermore, $D(G')$ has one crossing less than $D(G)$, $D(G')$ has the desired number of almost full
crossings,
and \textit{(i)} and \textit{(ii)} immediately follow.

 Assume that $S$ is a minimal cut of $G'$ and $v\in S$. Then there are $u,w\in N_{G'}(v)$ such that $u$ and $w$ belong to distinct components of $G'\setminus S$, thus,  $G'[N_{G'}(v)\setminus S]$ has to be disconnected. Since the neighborhood $N_{G'}(v)$ of $v$ in $G'$ forms an induced cycle (if $X$ is full) or an induced path (if $X$ is almost full) on four vertices (note that $G$ is simple), $S\cap
N_{G'}(v)\neq\emptyset$.

If $G'$ is not 3-connected, then $G'$ has a minimal cut $S$ such that $|S|\leq 2$ and $v\in S$. It follows that $X$ is not full and the subgraph of $G'$ spanned by $N_{G'}(v)$ is a path $P$ with one of its inner vertices in $S$. But then
both inner vertices of $P$ form a $2$-cut of $G$, in contradiction to the $3$-connectedness of $G$ and \textit{(iii)} is proved.

Now, we prove \textit{(iv)}.
First, let $S$ be a $3$-cut of $G'$ with $S\subseteq V(G)$. 
Then there are components $H_1$ and $H_2$ of $G'\setminus S$ with $v\in V(H_1)$. Because $N_{G'}(v)\subseteq V(H_1)\cup S$ and $|S|=3$, at most three of the four neighbors of $v$ belong to $S$. Hence, $S$ is a $3$-cut of $G'\setminus \{v\}$ and also of $G$. 
Next, assume 
that $G'$ contains a  $3$-cut $S$ with $v\in S$.
Let $e=xy$ and $e'=x'y'$ be two edges of the chosen crossing of $D(G)$, i.\,e.\ $N_{G'}(v)=\{x,y,x',y' \}$.

Case~1: \textit{$G'[\{x,y,x',y' \}]$ is a cycle on $4$ vertices. }\\
Because $G'[\{x,y,x',y' \}\setminus S]$ is disconnected, it follows that $S$ contains two independent neighbors of $v$, say
$S=\{v,x',y'\}$.
Thus, $G'-S$ has two components each containing a vertex of $\{x,y\}$.
If $G'-S$ has a further component $H$, then $V(H)\cap\{x,y,x',y' \}=\emptyset$ and $\{x',y'\}$ is a 2-cut of $G'$, a contradiction.
It is easy to see that there is  either an open Jordan curve  $J$  of the plane connecting $x'$ and $y'$ such that $J\cap
D(G)=\{x',y'\}$
or two edges, one from each component of $G'-S$, cross each other. The latter case cannot occur since the vertices of two crossing edges are connected in $G'$. Thus, 
if the edge $x'y'$ is replaced with $J$, then we get a drawing $D'(G)$  of $G$ with fewer crossings than $D(G)$, a
contradiction to the choice of $D(G)$.
It follows that Case~1 does not occur.

Case~2: \textit{$G'[\{x,y,x',y' \}]$ is a path on $4$ vertices. }\\
Without loss of generality assume $yy'\notin E(G)$. Because $G'[\{x,y,x',y' \}\setminus S]$ is disconnected, it follows that $S$ contains $x$ or $x'$.
Because of symmetry, let $x\in S$, i.\,e.\ $S=\{v,x,z\}$ with $z\in V(G)$.

Case~2.1: \textit{$z\notin\{x,y,x',y'\}$.}\\
With a similar argument as in Case~1, $G'-S$ has exactly two components $H_1$ and $H_2$ with $x',y\in V(H_1)$ and $y'\in V(H_2)$. 
It is easy to see that $S'=\{x,x',z\}$ is a $3$-cut of $G$. Moreover, $G-S'$ has two components $H_1-x'$ and $H_2$ each containing one vertex from $\{y,y'\}$.

Case~2.2: \textit{$z=y$.}\\
Then $S=\{v,x,y\}$ and we use the same arguments as in Case~1 for a contradiction.

Case~2.3: \textit{$z=x'$.}\\
Then $S=\{v,x,x'\}$ and, since $G'-S$ is disconnected, $G-\{x,x'\}$ is disconnected, contradicting the 3-connectedness of $G$. \hfill $\Box$

\bigskip

\begin{lemma}\label{lem:ind}
Let $G$ be a $1$-planar graph, $D(G)$ a $1$-planar drawing of $G$, and $e=xy$ and $e'=x'y'$ two crossing edges of $D(G)$.
Moreover, let $G'$ be obtained from $G$ by turning the crossing of  $e$ and $e'$ into a new $4$-valent vertex $v$.

a) If this crossing is full and $G'$ has a hamiltonian cycle $C'$ (hamiltonian path $P'$), then $G$ is hamiltonian (traceable).

b) If this crossing is almost full with $xx'\notin E(G)$, and $G'$ has a hamiltonian cycle $C'$ (hamiltonian path $P'$) not containing both $vx$ and $vx'$, then $G$ is
hamiltonian (traceable).
\end{lemma}

{\bf Proof of \thref{lem:ind}.}\\
Let $vu$ and $vw$ be adjacent edges of $C'$. In either case we have $uw\in E(G)$, so replacing the subpath $uvw$ of $C'$ with the edge $uw$ leads to a hamiltonian cycle of $G$. The same arguments hold for $P'$ if $v$ is not an end vertex of $P'$. If it is, simply remove it and we obtain the desired hamiltonian path in $G$.
\hfill $\Box$

\bigskip

The following lemma is a result of Brinkmann and the last author~\cite{bz}:

\begin{lemma}[Theorem~16 and Corollary~17 from~\cite{bz}]\label{lem:gunnarcarol}
Each  $3$-connected planar graph with at most three $3$-cuts is hamiltonian and each  $3$-connected planar graph with at most four $3$-cuts is traceable. 
\end{lemma}

\bigskip

{\bf Proof of}
\ifrestated
\vspace{-1.6ex}

\thmthree*
\else
{\bf\thref{thm:3cut}. }
\fi

Let $G_1$ be a $3$-connected locally maximal $1$-planar graph with at most three $3$-cuts.
We define a sequence of locally maximal $1$-planar graphs $G_1, G_2, \dots$, where for all $i\ge 1$, $G_{i+1}$ is the graph $G'$ if $G_i$ is the
graph $G$ according to \thref{lem:struct} (with $t=0$).
By \thref{lem:struct}, there is an index $k$ such that $G_k$ is planar and $3$-connected with at most three $3$-cuts;  no further $3$-cut appears since all crossings of $G_1$ are full.
By \thref{lem:gunnarcarol}, $G_k$ is hamiltonian.
Applying assertion a) of \thref{lem:ind} repeatedly implies that $G_1$ is hamiltonian.

In the same spirit, let now $G_1$ be a $3$-connected locally maximal $1$-planar graph with at most four $3$-cuts. Define a sequence $G_1, G_2, \dots$ as above. By \thref{lem:struct}, there is an index $k$ such that $G_k$ is planar and $3$-connected with at most four $3$-cuts. 
By \thref{lem:gunnarcarol}, $G_k$ is traceable. Again we apply assertion~a) of \thref{lem:ind} repeatedly and obtain that $G_1$ is traceable.
 \hfill $\blacksquare$

\bigskip

{\bf Proof of}
\ifrestated
\vspace{-1.6ex}

\thmfour*
\else
{\bf\thref{thm:almost}. }
\fi

Let $G_1$ be a $4$-connected $1$-planar graph which has a weakly locally maximal $1$-planar drawing with at most three almost full crossings.
Among all weakly locally maximal $1$-planar drawings of $G_1$ with at most three almost full crossings, let $D(G_1)$ be chosen with minimum number
of crossings.
If the number $s$ of almost full crossings in $D(G_1)$ is zero, then $G_1$ is hamiltonian by \thref{thm:3cut}.

We assume $s\ge 1$, consider an almost full crossing $X$ of $D(G_1)$ and apply \thref{lem:struct} to this crossing with $G=G_1$ and $t=s$, and
obtain $G_1'=G'$ with the new added vertex $v_1=v$. Obviously, $G_1'$ has a drawing with at most $s-1$ almost full crossings.
Since $G_1$ is $4$-connected, $G_1'$ is $4$-connected 
by  \thref{lem:struct}. 
Let $G_2$ be obtained from $G_1'$ by adding  a vertex $u_1$, the edge $u_1v_1$ and the two
edges connecting $u_1$ with both 2-valent vertices of the path  $G_1'[N_{G_1'}(v_1)]$ (see \thref{lem:struct}  and Figure~\ref{pic:u}).
Then, $G_2$ is $3$-connected and $N_{G_2}(u_1)$ is the only 3-cut of $G_2$.
Furthermore, $G_2$ has a weakly locally maximal $1$-planar drawing with $s-1$ almost full crossings.

\begin{figure}[h]
\begin{minipage}{.5\linewidth}
\begin{center}
\begin{tikzpicture}[scale=0.4]

\fill (-2,2) circle (0.2cm);
\fill (2,2) circle (0.2cm);
\fill (-2,-2) circle (0.2cm);
\fill (2,-2) circle (0.2cm);

\draw[thick] (-2,2) -- (2,2);
\draw[thick] (-2,2) -- (-2,-2);
\draw[thick] (-2,2) -- (2,-2);
\draw[thick] (-2,-2) -- (2,2);
\draw[thick] (-2,-2) -- (2,-2);

\draw (-2,3.4) node {$b\vphantom{y'}$};
\draw (2,3.4) node {$a\vphantom{y'}$};
\draw (-2,-3.4) node {$c\vphantom{y'}$};
\draw (2,-3.4) node {$d\vphantom{y'}$};

\fill (0,0) circle (0.2cm);
\draw (1,0) node {$v_i$};

\fill (-1,0) circle (0.2cm);
\draw (-2.6,0) node {$u_i$};

\draw[thick,dashed] (-2,2) -- (-1,0);
\draw[thick,dashed] (-2,-2) -- (-1,0);
\draw[thick,dashed] (0,0) -- (-1,0);

\end{tikzpicture}
\end{center}

\caption{}\label{pic:u}
\end{minipage}
\begin{minipage}{.48\linewidth}
\begin{center}
\begin{tikzpicture}[scale=0.4]

\fill (-2,2) circle (0.2cm);
\fill (2,2) circle (0.2cm);
\fill (-2,-2) circle (0.2cm);
\fill (2,-2) circle (0.2cm);

\fill (-6,2) circle (0.2cm);
\fill (-6,-2) circle (0.2cm);

\draw (-2,3.4) node {$b_i=b=b_j\vphantom{y'}$};
\draw (-2,-3.4) node {$c_i=c=c_j\vphantom{y'}$};

\draw (-6.9,2) node {$a_i\vphantom{y'}$};
\draw (-6.9,-2) node {$d_i\vphantom{y'}$};
\draw (2.9,2) node {$a_j\vphantom{y'}$};
\draw (2.9,-2) node {$d_j\vphantom{y'}$};

\draw[thick] (-2,2) -- (2,2);
\draw[thick] (-2,2) -- (-2,-2);
\draw[thick] (-2,2) -- (2,-2);
\draw[thick] (-2,-2) -- (2,2);
\draw[thick] (-2,-2) -- (2,-2);

\draw[thick] (-2,2) -- (-6,2);
\draw[thick] (-2,2) -- (-6,-2);
\draw[thick] (-2,-2) -- (-6,2);
\draw[thick] (-2,-2) -- (-6,-2);

\fill (0,0) circle (0.2cm);
\draw (1,0) node {$v_j$};

\fill (-4,0) circle (0.2cm);
\draw (-5,0) node {$v_i$};

\fill (-1,0) circle (0.2cm);
\fill (-3,0) circle (0.2cm);

\draw[ultra thick,line width=2.2pt] (-2,2) -- (-1,0);
\draw[ultra thick,line width=2.2pt] (-2,-2) -- (-1,0);
\draw[thick,] (0,0) -- (-1,0);

\draw[ultra thick,line width=2.2pt] (-2,2) -- (-3,0);
\draw[ultra thick,line width=2.2pt] (-2,-2) -- (-3,0);
\draw[thick,] (-4,0) -- (-3,0);

\draw[ultra thick,line width=2.2pt] (0,0) -- (2,2);
\draw[ultra thick,line width=2.2pt] (0,0) -- (2,-2);
\draw[ultra thick,line width=2.2pt] (-4,0) -- (-6,2);
\draw[ultra thick,line width=2.2pt] (-4,0) -- (-6,-2);

\end{tikzpicture}
\end{center}

\caption{}\label{pic:double_u}
\end{minipage}
\end{figure}

Note that a hamiltonian cycle of $G_2$ (if it exists) leads to a hamiltonian cycle of $G_1'=G_2\setminus \{u_1\}$ containing at least one edge of $G_1'[N_{G_2}(u_1)]$.

If $s=1$, then let $H=G_2$.
Otherwise, we repeat this step $s-1$ times and obtain a graph $H=G_3$ or $H=G_4$.
$H$ is  $3$-connected locally maximal $1$-planar. 
Assume there is a $3$-cut $S$ in $G_{i+1}$ which is not a $3$-cut of $G_i$ for $i\in \{1,\dots,s\}$. Then $S=N_{G_{i+1}}(u_i)$ or $S\neq N_{G_{i+1}}(u_i)$ and $v_i\in S$. 
In the second case, by \thref{lem:struct}, there is a $3$-cut in $G_i$ separating two vertices of $N_{G_{i+1}}(v_i)\setminus\{u_i\}$, a contradiction. 
Hence, $H$ has exactly $s$ $3$-cuts, namely the neighborhoods of $u_1,\dots, u_s$.
Since $s\le 3$, $H$ is hamiltonian by \thref{thm:3cut} and $G=H-\{u_1,\dots, u_s\}$ is hamiltonian because the neighborhoods of $u_1,\dots,
u_s$ are complete in $H$.
By the previous remark, we may assume that $G$ contains a hamiltonian cycle $C$ such that
\begin{align}\label{star}
E(G[N_H(u_i)])\cap E(C)\neq \emptyset\text{ for $i=1,\dots,s$.} \tag{$\ast$}
\end{align}

Consider an arbitrary vertex $v\in\{ v_1,\dots, v_s\}$ and let $\{a,b,c,d\}$ be the vertex set of the induced path on $N_G(v)$ in this order.

If $av,dv\in E(C)$, then $bc\in E(C)$ by property (\ref{star}).
In this case let $P_1$ and $P_2$ be the subpaths of $C$ obtained by removing $v$ and the edge $bc$ from $C$.
If $P_1$ connects $ab$ and $P_2$ connects $cd$, then the cycle obtained from $P_1$, $P_2$, $ac$, and $bd$ is a hamiltonian cycle of the graph
$G_v$ obtained from $G$ by deleting $v$ and adding the edges $ac$ and $bd$.
If $P_1$ connects $ac$ and $P_2$ connects $bd$, then the cycle obtained from $P_1$, $P_2$, $ab$, and $cd$ is a hamiltonian cycle of the graph
$G_v$.
If not both edges $av,dv$ belong to $C$, then $G_v$ is hamiltonian by assertion b) of \thref{lem:ind}.

Repeating this step $s$ times, we get rid of $v_1,\dots, v_s$, the resulting graph is $G_1$ and the existence of a hamiltonian cycle of $G_1$ is
shown. Note that if there exist distinct $v_i$ and $v_j$ sharing the same $bc$, 
 then $C$ misses at least one edge of $a_iv_i, d_iv_i, a_jv_j, d_jv_j$ (see Figure~\ref{pic:double_u}), since otherwise the edges $bu_i,u_ic,cu_j,u_jb$ of $C$ form a cycle. 

The proof that any $4$-connected $1$-planar graph which has a weakly locally maximal $1$-planar drawing with at most four almost full crossings is traceable uses very similar arguments and is therefore omitted. \hfill $\blacksquare$

\bigskip

{\bf Proof of}
\ifrestated
\vspace{-1.6ex}

\thmfive*
\else
{\bf\thref{Subgraph}. }
\fi

Among all locally maximal $1$-planar drawings  of $G$ let $D(G)$ be chosen such that the number of crossings in $D(G)$ is minimal.
If two edges of $D(G)$ cross each other, then remove an arbitrary one of them and let $H$ be the resulting graph.
Obviously, $H$ is plane and a spanning subgraph of $D(G)$ (and of $G$).

It remains to show that $H$ is $3$-connected.
Assume $H$ is not $3$-connected and, therefore, let $S\subset V(H)$ be a cut of $H$ with $|S|\le 2$ such that $H_1,\dots,H_k$ ($k\ge 2$)
are the components of $H-S$.
Since $D(G)-S$ is connected, there are at least $k-1$ {\em connecting} edges $xy\in E(D(G))\setminus E(H)$ with $x\in V(H_i)$ and $y\in V(H_j)$
for suitable $i,j\in \{1,\dots,k\}$ with $i\neq j$.
 These edges are crossed in $D(G)$ by edges from $E(H)$.

Let $x'y'$ be the edge crossing some connecting edge $xy\in E(D(G))\setminus E(H)$ in $D(G)$.

Since $D(G)$ is locally maximal, $\{x,y,x',y'\}$ induces a complete subgraph of $G$, thus, both $x'$ and $y'$ are common neighbors of $x$ and
$y$.

It follows that $S=\{x',y'\}$. Hence $xy$ is the only connecting edge (another connecting edge  would also cross $x'y'$) and therefore $k=2$.
We argue as in Case 1 of the proof of  \thref{lem:struct}: there is an open Jordan curve  $J$  of the plane connecting $x'$ and $y'$ such that $J\cap
D(G)=\{x',y'\}$ and, if the edge $x'y'$ is replaced with $J$, then we get a drawing $D'(G)$  with fewer crossings than $D(G)$, a
contradiction.
 \hfill $\blacksquare$

\section{Overview}

We end this paper with a tabular overview of hamiltonian properties of various families of graphs that we have discussed.

\begin{table}[h]
\definecolor{yes}{rgb}{0.5,1.0,0.5}
\definecolor{no}{rgb}{1.0,0.5,0.5}
\newcommand{\yes}{\cellcolor{yes} \ding{51}}
\newcommand{\no}{\cellcolor{no} \ding{55}}
\footnotesize
\begin{center}
\begin{tabular}{cccccccc}
\toprule
 & Maximal & Planar & Optimal  & Maximal  & Locally  & Weakly  & 1-planar\\
 & planar  &        & 1-planar & 1-planar & maximal  & locally & \\
 &         &        &          &          & 1-planar & maximal & \\
 &         &        &          &          &          & 1-planar & \\
\midrule
3 & \no & \no & \rule[0.8mm]{3.5mm}{0.4mm} & \no\makebox[0pt][l]{(D1)} & \no & \no & \no \\
4 & \yes\makebox[0pt][l]{(A)} & \yes\makebox[0pt][l]{(B)} & \yes\makebox[0pt][l]{(C)} & \yes\makebox[0pt][l]{(D3)} & \yes\makebox[0pt][l]{\hspace{-0.3mm}(D3)} & \no  & \no\makebox[0pt][l]{(C)} \\
5 & \yes & \yes & \yes & \yes &  \yes & \no\makebox[0pt][l]{(D2)} & \no\makebox[0pt][l]{\hspace{-0.3mm}(D2)} \\
6 & \rule[0.8mm]{3.5mm}{0.4mm} & \rule[0.8mm]{3.5mm}{0.4mm} &  \yes & \yes & \yes & \textbf{?} & \textbf{?} \\
7 & \rule[0.8mm]{3.5mm}{0.4mm} & \rule[0.8mm]{3.5mm}{0.4mm} &  \yes & \yes\makebox[0pt][l]{(C)} & \yes & \textbf{?} & \textbf{?} \\
\bottomrule
\end{tabular}
\medskip

\begin{tabular}{m{-2mm} l m{2mm} l}
A. & Whitney \cite{whitney1931theorem} & D1. & This paper, Theorem 1\\
B. & Tutte \cite{tutte1956theorem} & D2. & This paper, Theorem 2\\
C. & Hud{\'a}k, Madaras, Suzuki \cite{hudak2012properties} \qquad & D3. & This paper, Theorems 3 and 4\\

\end{tabular}
\end{center}
\caption{Hamiltonicity of planar and 1-planar graphs, as well as some related families, listed by connectedness ranging from 3 to 7 (the maximum admissible value for $1$-planar graphs). Green cells (marked \ding{51}) indicate that every graph with the specified connectedness is hamiltonian, red cells (marked \ding{55}) signify that there exist such graphs which are not hamiltonian, question marks designate open problems, and ``\rule[0.8mm]{3.5mm}{0.4mm}'' stands for an impossible combination of properties.} \label{tab:overview}
\end{table}

{\bf Acknowledgment.} 
We thank the referees, whose constructive suggestions helped to improve the quality of this article.


\begin{thebibliography}{1}

\bibitem{Brandenburg1}
C. Bachmaier, F.J. Brandenburg, K. Hanauer, D. Neuwirth, and J. Reislhuber,
NIC-planar graphs, 
{\em Discrete Appl. Math.} 232(2017)23--40.

\bibitem{Barnette}
D.W. Barnette,
Trees in polyhedral graphs,
{\em Canad. J. Math.} 18(1966)731--736.

\bibitem{Brandenburg}
F.J. Brandenburg, D. Eppstein, A. Glei\ss ner, M.T. Goodrich, K. Hanauer, and J. Reislhuber,
On the density of maximal 1-planar graphs,
In: GD 2012 (eds.: M. van Kreveld and B. Speckmann), LNCS 7704, pp.~327--338, Springer, 2013.

\bibitem{bz}
G. Brinkmann and C.T. Zamfirescu,
Polyhedra with few 3-cuts are hamiltonian,
{\em Electron. J. Combin.} 26.1(2019).

\bibitem{ChenYu}
G. Chen and X. Yu,
Long cycles in 3-connected graphs,
{\em J. Combin. Theory, Ser. B} 86(2002)80--99.

\bibitem{fabrici2007structure}
I. Fabrici and T. Madaras,
The structure of 1-planar graphs,
{\em Discrete Math.} 307(2007)854--865.

\bibitem{Fleischner}
H. Fleischner, The uniquely embeddable planar graphs, 
{\em Discrete Math.} 4(1973)347--358.

\bibitem{Gao}
Z. Gao,
2-connected coverings of bounded degree in 3-connected graphs,
{\em J. Graph Theory} 20(1995)327--338.

\bibitem{GW}
B. Gr\"unbaum and H. Walther,
Shortness exponents of families of graphs,
{\em J. Combin. Theory, Ser. A} 14(1973)364--385.

\bibitem{hudak2012properties}
D. Hud{\'a}k, T. Madaras, and Y. Suzuki,
On properties of maximal 1-planar graphs,
{\em Discuss. Math. Graph Theory} 32(2012)737--747.

\bibitem{Korzhik}
V.P. Korzhik,
Minimal non-$1$-planar graphs,
{\em  Discrete Math.} 308(2008)1319--1327.

\bibitem{KorzhikMohar}
V.P. Korzhik and B. Mohar,
Minimal obstructions for 1-immersions and hardness of 1-planarity testing,
{\em J. Graph Theory} 72(2013)30--71.

\bibitem{MoonMoser}
J.W. Moon and L. Moser,
Simple paths on polyhedra,
{\em Pacific J. Math.} 13(1963)629--631.

\bibitem{pach1997graphs}
J. Pach and G. T{\'o}th,
Graphs drawn with few crossings per edge,
{\em Combinatorica} 17(1997)427--439.

\bibitem{ringel1965sechsfarbenproblem}
G. Ringel,
Ein Sechsfarbenproblem auf der Kugel,
{\em Abh. Math. Sem. Univ. Hamburg} 29(1965)107--117.

\bibitem{tutte1956theorem}
W.T. Tutte,
A theorem on planar graphs,
{\em Trans. Amer. Math. Soc.} 82(1956)99--116.

\bibitem{whitney1931theorem}
H. Whitney,
A Theorem on Graphs,
{\em Ann. Math.} 32(1931)378--390.

\bibitem{whitney1932theorem}
H. Whitney,
Congruent graphs and the connectivity of graphs, 
{\em Am. J. Math.} 54(1932)150--168.

\end{thebibliography}
\end{document}